\documentclass{llncs}

\usepackage[utf8x]{inputenc}
\usepackage{verbatim}
\usepackage{float}

\usepackage{amsmath}

\usepackage{amsthm}
\usepackage{amssymb}
\usepackage{algpseudocode}
\usepackage{array}
\usepackage{graphicx}
\usepackage{longtable}
\usepackage{enumerate}
\usepackage{siunitx}
\usepackage{bookmark}
\usepackage{hyperref}
\hypersetup{
  pdftitle = {Integer Complexity: Experimental and Analytical Results},
  pdfauthor = {Jānis Iraids, Kaspars Balodis, Juris Čerņenoks, Mārtiņš Opmanis, Rihards Opmanis, Kārlis Podnieks},
  pdfkeywords = {integer complexity, shortest expression, experimental mathematics},
  pdfstartview=FitH,
  unicode=true
}

\begin{document}

\title{Integer Complexity: Experimental and Analytical Results}
\author{\texorpdfstring{Jānis Iraids\inst{1} \and Kaspars Balodis\inst{2} \and Juris Čerņenoks\inst{2} \and Mārtiņš Opmanis\inst{1} \and Rihards Opmanis\inst{1} \and Kārlis Podnieks\inst{1}}{Jānis Iraids, Kaspars Balodis, Juris Čerņenoks, Mārtiņš Opmanis, Rihards Opmanis, Kārlis Podnieks}}
\titlerunning{Integer Complexity}
\institute{Institute of Mathematics and Computer Science, University of Latvia, Raiņa bulvāris 29, Riga, LV-1459, Latvia \and University of Latvia, Raiņa bulvāris 19,  Riga, LV-1586, Latvia}

\makeatletter
\providecommand*{\toclevel@title}{0}
\providecommand*{\toclevel@author}{0}
\makeatother

\theoremstyle{plain}
\newtheorem{hypo}{Hypothesis}
\newtheorem{obs}{Observation}

\newcommand{\f}[1]{ \left\| #1 \right\| }
\newcommand{\fl}[1]{ \f{ #1 }_{\log} }

\maketitle

\bookmarksetup{startatroot}

\begin{abstract}
We consider representing of natural numbers by arithmetical expressions using ones, addition, multiplication and parentheses. The (integer) complexity of $n$ -- denoted by $\f{n}$ -- is defined as the number of ones in the shortest expressions representing $n$. We arrive here very soon at the problems that are easy to formulate, but (it seems) extremely hard to solve. In this paper we represent our attempts to explore the field by means of experimental mathematics. Having computed the values of $\f{n}$ up to $10^{12}$ we present our observations. One of them (if true) implies that there is an infinite number of Sophie Germain primes, and even that there is an infinite number of Cunningham chains of length 4 (at least). We prove also some analytical results about integer complexity.
\end{abstract}

\section{Introduction}
\pdfbookmark[1]{1. Introduction}{Introduction}

The field explored in this paper is represented most famously in F26 of Guy \cite{f26}, and as the sequence A005245 in “The On-Line Encyclopedia of Integer Sequences”\cite{oeisA005245}.

We consider representing of natural numbers by arithmetical expressions using ones, addition, multiplication and parentheses. Let's call this “representing numbers in basis $\{1,+,\cdot\}$”. For example,
{
\allowdisplaybreaks
\begin{eqnarray*}
2&=&1+1;\\
3&=&1+1+1;\\
4&=&1+1+1+1=2\cdot2=(1+1)\cdot (1+1);\\
5&=&1+1+1+1+1=1+2\cdot 2=1+(1+1)\cdot (1+1);\\
6&=&1+1+1+1+1+1=2\cdot 3=(1+1)\cdot (1+1+1);\\
7&=&1+1+1+1+1+1+1=1+2\cdot 3=1+(1+1)\cdot (1+1+1);\\
8&=&1+1+1+1+1+1+1+1=\\
 &=&2+2\cdot 3=1+1+(1+1)\cdot (1+1+1)=\\
 &=&2\cdot 2\cdot 2=(1+1)\cdot (1+1)\cdot (1+1);\\
9&=&1+1+1+1+1+1+1+1+1=3+2\cdot 3=1+2\cdot 2\cdot 2=\\
 &=&3\cdot 3=(1+1+1)\cdot (1+1+1);\\
10&=&1+1+\ldots +1=2\cdot 2+2\cdot 3=2+2\cdot 2\cdot 2=\\
&=&2\cdot 5=(1+1)\cdot (1+1+1+1+1);\\
11&=&\ldots =2+3\cdot 3=1+2\cdot 5;\\
12&=&\ldots =2\cdot 2\cdot 3;\\
13&=&\ldots =2\cdot 2+3\cdot 3=3+2\cdot 5;\\
14&=&\ldots =2\cdot (1+2\cdot 3);\\
\ldots
\end{eqnarray*}
}

As we see, most numbers can be represented in several ways that may differ in size. For example, the number 8 is represented above by three different expressions containing 8, 7 and 6 ones respectively.

We will measure the size of an expression by the number of ones it contains. We don't need counting neither of operations (if an expression contains $k$ ones, then it contains $k–1$ operations), nor of parentheses (the postfix notation might be used).

The size of the shortest expressions representing a particular number $n$ can be considered as the “complexity” of $n$. Hence, the term “integer complexity”. Some numbers allow for several shortest expressions (for examples, see above: 4, 5, 10, 11, 13). 

\begin{definition}
Let's denote by $\f{n}$ the number of ones in the shortest expressions representing $n$ in basis $\{1,+,\cdot \}$. We will call it the \textbf{integer complexity} of $n$.
\end{definition}

For example, as we see above:

\[
\begin{array}{llll}
\f{1}=1; &\qquad \f{5}=5; &\qquad \f{9}=6;  &\qquad \f{13}=8; \\
\f{2}=2; &\qquad \f{6}=5; &\qquad \f{10}=7; &\qquad \f{14}=8; \\
\f{3}=3; &\qquad \f{7}=6; &\qquad \f{11}=8; &\qquad \f{15}=8; \\
\f{4}=4; &\qquad \f{8}=6; &\qquad \f{12}=7; &\qquad \ldots
\end{array}
\]

This definition corresponds to the sequence A005245 in “The On-Line Encyclopedia of Integer Sequences”. In \cite{f26}, $\f{n}$ is denoted by $f(n)$, the notation $\f{n}$ is due to Arias de Reyna \cite{arias}.

In a similar fashion, representation of natural numbers in other bases, for example, $\{1,+,\cdot, -\}$, $\{1,+,\cdot ,\uparrow \}$ and $\{1, +, \cdot, –, \uparrow \}$ could also be considered (sequences A091333, A025280 and A091334 \cite{oeisA091333,oeisA025280,oeisA091334}, $\uparrow $ stands for exponentiation).

As a function of $n$, in average, $\f{n}$ is growing logarithmically, namely, one can prove easily:
\begin{theorem}
\label{cbounds}
For all $n>1$, \[3\log_3 n \leq \f{n} \leq 3\log_2 n \approx 4.755 \log_3 n\]
\end{theorem}
In \cite{f26}, Guy attributes this result to Dan Coppersmith.

The lower bound of Theorem \ref{cbounds} is reached by infinitely many numbers $n$, exactly – by the powers of three. For example,
\[3^4 = 81 = (1+1+1)\cdot (1+1+1)\cdot (1+1+1)\cdot (1+1+1),\]
and, indeed, products of $1+1+1$s is the best way of representing powers of three in basis $\{1,+,\cdot \}$:
\begin{theorem}
\label{cbounds2}
$\f{n}= 3\log_3n$, if and only if $n=3^b$ for some $b\geq 1$. In particular, for all $b\geq 1$, $\f{3^b}=3b$ (moreover, the product of $1+1+1$s is shorter than any other representation of $3^b$).
\end{theorem}

Similarly, the product of $1+1$s seems to be the best way of representing powers of two in basis $\{1,+,\cdot \}$. For example,
\[2^5=32=(1+1)\cdot (1+1)\cdot (1+1)\cdot (1+1)\cdot (1+1).\]

\begin{hypo}
\label{powtwohypo}
For all $a\geq 1$, $\f{2^a}=2a$ (moreover, the product of $1+1$s is shorter than any other representation of $2^a$).
\end{hypo}

$\f{2^a}=2a$ is true for all powers $2^a<10^{12}$, i.e. for all $a$, $0<a\leq 39$ – as verified by Jānis Iraids. We consider proving or disproving of Hypothesis \ref{powtwohypo} as one of the biggest challenges of number theory.

The upper bound of Theorem \ref{cbounds} doesn't seem to be exact. As observed by Rawsthorne \cite{raw}, the “worst” value of $\frac{\f{n}}{\log_3 n}$ seems to be $\approx 3.928$ (and not $4.755$ of Theorem \ref{cbounds}), and it is reached only by a single number, namely, by $1439$ ($\f{1439}=26$):
\begin{hypo}
\label{maxeffhypo}
For all $n>1$,
\[\f{n} \leq \frac{\f{1439}}{\log_3 1439} \log_3 n \approx 3.928\log_3 n\]
\end{hypo}
Hypothesis \ref{maxeffhypo} is true for all $n\leq 10^{12}$ – as verified by Jānis Iraids.

Thus, we arrive here very soon at the problems that are easy to formulate, but (it seems) extremely hard to solve. In this paper we represent our attempts to explore the field by means of experimental mathematics. We managed to prove analytically only a few of the experimental observations.

In Section 2 we explain the basic concepts and their simplest properties. Section 3 represents our analytical results. Section 4 considers algorithms allowing to calculate the values of $\f{n}$. The best of them was used by Jānis Iraids to calculate $\f{n}$ up to $n=10^{12}$. In Section 5 we present our experimental observations (several confirmed and refuted hypotheses included).

Wolfram|Alpha (by Wolfram Alpha LLC, 2009) was used for some of the computations in this paper.

\section{Basic concepts and related work}
\pdfbookmark[1]{2. Basic concepts and related work}{Basic concepts and related work}

\subsection{The largest and smallest numbers of complexity $n$}
\pdfbookmark[2]{2.1. The largest and smallest numbers of complexity n}{The largest and smallest numbers of complexity n}

For a given $n$, there exists only a finite number of expressions of size $\leq n$. Hence,

\begin{definition}
\label{e}
\begin{enumerate}
\item[a)] Let's denote by $E(n)$ the  \textbf{largest} $m$ such that $\f{m}=n$.
\item[b)] Let's denote by $E_k(n)$ the \textbf{ $k$-th largest} $m$ such that $\f{m}\leq n$ (if it exists). Thus, $E(n)=E_0(n)$.
\item[c)] Let's denote by $e(n)$ the  \textbf{smallest} $m$ such that $\f{m}=n$.
\end{enumerate}
\end{definition}

In this definition, $E(n)$ corresponds to the sequence A000792\cite{oeisA000792}, and $e(n)$ – to A005520\cite{oeisA005520}.

\begin{proposition}
\label{emono}
As a function of $n$, $E(n)$ is monotonically increasing.
\end{proposition}

\begin{proof}
Take one of the shortest expressions for $E(x)$ and append $+1$ to it. Now, it contains $x+1$ ones, and its value is $E(x)+1$. Therefore the greatest number of complexity $x+1$ is no less than $E(x)+1$.
\end{proof}

\begin{proposition}
As a function of $n$, $e(n)$ is monotonically increasing.
\end{proposition}

\begin{proof}
Assume that there exists $n$, such that $e(n)<e(n-1)$. Take the smallest such $n$. Then it would mean that before $e(n)$ there are no numbers of complexity $\geq n$, in particular, $\f{e(n)}-\f{e(n)-1}>1$. But it is impossible since $e(n)$ can be written as $e(n)+1$.
\end{proof}

\begin{proposition}
\label{xleqefx}
For all $n\geq 1$, $n \leq E(\f{n})$.
\end{proposition}

\begin{proof}
By definition -- $E(y)$ is the greatest number whose complexity is $y$. Thus for all $n$, if $\f{n}=y$ then $n \leq E(y)$.
\end{proof}

\begin{proposition}
\label{exeyleqexy}
For all $x, y\geq 1$, $E(x)\cdot E(y) \leq E(x+y)$.
\end{proposition}

\begin{proof}
Take the product of shortest expressions for $E(x)$ and $E(y)$. The value of this product is $E(x)\cdot E(y)$ and it contains $x+y$ ones, so it cannot be greater than $E(x+y)$.
\end{proof}

\begin{theorem}
\label{eexp}
For all $k\geq 0$:
\[E(3k+2) = 2\cdot 3^k;\]
\[E(3k+3) = 3\cdot 3^k;\]
\[E(3k+4) = 4\cdot 3^k.\]
\end{theorem}

Guy \cite{guy} attributes this result to John L. Selfridge.

\begin{theorem}
For all $n\geq 8$: \[E_2(n) = \frac{8}{9}E(n).\]
\end{theorem}

This result is due to Rawsthorne\cite{raw}.

The behaviour of $e(n)$ appears to be more complicated, for details see Section \ref{subsec:eexp}.

\begin{lemma}
If $a+b=n$ and $\f{a}+\f{b}=\f{n}$ and $\f{a}\leq \f{b}$ and given $\f{n}\leq N$, then $E(\f{a}) \leq \frac{n-\sqrt{n^2-4E(N)}}{2}$.
\end{lemma}

\begin{proof}
From $a+b=n$ we get $E(\f{a})+E(\f{b})\geq n$ using Proposition \ref{xleqefx}.
Further, as $\f{a}+\f{b}=\f{n}$, we get $E(\f{a})+\frac{E(\f{n})}{E(\f{a})}\geq n$ by Proposition \ref{exeyleqexy}. Still more due to the monotonicity of $E(x)$ we can substitute for the estimate of $n$'s complexity: $E(\f{a})+\frac{E(N)}{E(\f{a})}\geq n$. For convenience let us denote $E(\f{a})$ by $x$ and $E(N)$ by $y$ obtaining $x+\frac{y}{x}\geq n$. Solving the quadratic inequality for $x$ we get $x \leq \frac{n-\sqrt{n^2-4y}}{2}$. To complete the proof insert back the original values of $x$ and $y$.
\end{proof}

\begin{corollary}
\label{smallestaddend}
For $n\geq 29$, if $a+b=n$ and $\f{a}+\f{b}=\f{n}$ and $\f{a}\leq \f{b}$, then $a\leq 2n^{\log_2{3}-1} \approx 2n^{0.585}$.
\end{corollary}

\begin{proof}
By Theorem \ref{cbounds}, $N \leq 3\log_2{n}$. Furthermore, we use the convenient fact that $E(n) \leq 3^{\frac{n}{3}}$.
\[\begin{split} a &\leq E(\f{a}) \leq \frac{n-\sqrt{n^2-4E(N)}}{2} \leq \\
&\leq \frac{n-\sqrt{n^2-4\cdot 3^{\log_2{n}}}}{2} \leq \frac{n(1-\sqrt{1-4\cdot n^{\log_2{3}-2}})}{2} \leq \\
&\leq \frac{n(1-(1-4\cdot n^{\log_2{3}-2}))}{2} \leq 2n^{\log_2{3}-1}. \end{split}\]
\end{proof}

A similar proof of Corollary \ref{smallestaddend} is given in \cite{vivek}.

\subsection{Ranking numbers}
\label{subsec:ranking}
\pdfbookmark[2]{2.2. Ranking numbers}{Ranking numbers}

Consider an expression in basis $\{1,+,\cdot \}$ drawn as a rooted $n$-ary tree, its leaves containing ones and inner nodes containing either $+$ or $\cdot $. Since both addition and multiplication is associative and commutative, let us merge any adjacent additions and multiplications. For example, the shortest expression for $5^6$ can be obtained as follows, it contains 29 ones (not $5 \cdot 6 = 30$, as one might expect):

\[5^6=15625=1+2^3 \cdot 3^2 \cdot 217=1+2^3 \cdot 3^2(1+2^3 \cdot 3^3)\]

\begin{figure}
  \centering
  \includegraphics[width=0.9\textwidth]{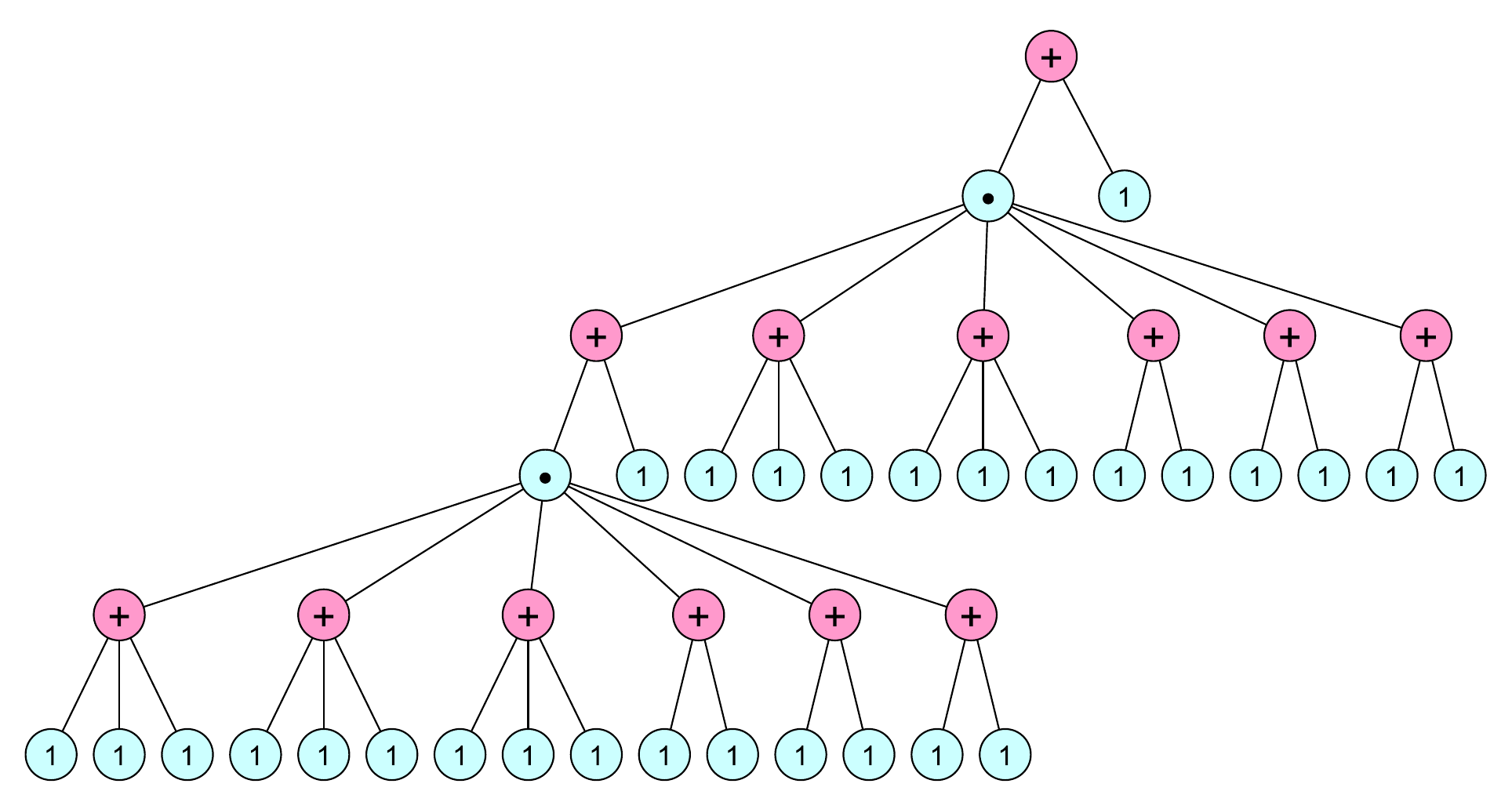}
  \caption{The tree of the shortest expression for $5^6$}
  \label{fig:rank}
\end{figure}

The corresponding tree is drawn in Figure \ref{fig:rank}.

\begin{definition}
The \textbf{height of an expression} is the height of the corresponding tree.
\end{definition}

In general, for a given $n$, there can be several shortest expressions of different height,  for example, of height 2 and 1:\[4=(1+1)\cdot (1+1) = 1+1+1+1,\] or of height 4, 3 and 2: \[10=(1+1)\cdot(1+(1+1)\cdot(1+1))=1+(1+1+1)\cdot(1+1+1)=(1+1)\cdot (1+1+1+1+1).\] We will be interested in classifying the shortest expressions by introducing expression height as another complexity measure of positive integers:

\begin{definition}
$rank(n)$ -- the \textbf{rank of a number} $n$ – is defined as the minimum height among the shortest expressions representing $n$.
\end{definition}

Examples:

\renewcommand{\labelitemi}{$\ast$}
\begin{itemize}
  \item the only number of rank $0$ is $1$;
  \item the only numbers of rank $1$ are: $2, 3, 4, 5$; 
  \item numbers of rank $2$: $6, 8, 9, 10, 12, 15, ...$; see Hypothesis \ref{genhypo} below;
  \item $rank(7)=rank(11)=rank(13)=3$;
  \item $rank(14)=4$.
\end{itemize}

For $n \leq 10^9$, $\max {rank(n)}=19$ and $\max \f{n}=67$. For other observations – see Section \ref{subsec:exprank}. 

Obviously $rank(n)=2$ for infinitely many numbers $n$. Namely, for all the powers of $3$.

\begin{hypo}
\label{strongunlrank}
For all $r>2$ there exists an infinite amount of numbers having rank $r$.
\end{hypo}

This hypothesis implies

\begin{hypo}
\label{unlrank}
As a function of $n$, $rank(n)$ is unlimited.
\end{hypo}

\begin{definition}
\label{rn}
Let's denote by $r(n)$ the \textbf{smallest} $m$ such that $rank(m)=n$.
\end{definition}

The values of $r(n)$ up to n=19 are represented in Table \ref{table:rn}.

Now, let's try exploring more closely the numbers of rank $2$.

\begin{lemma}
\label{primedec}
For $n>1$, the shortest expression of height $\leq 2$ representing $n$, is obtained from the prime decomposition $p_1\cdot \ldots \cdot p_k$ of $n$ as $(1+\ldots+1)\cdot (1+\ldots+1)\cdot \ldots\cdot (1+\ldots+1)$.
\end{lemma}

The number of ones in this expression, i.e. $p_1+\ldots +p_k$, is called, sometimes, the \textbf{integer logarithm} of $n$ (see also the sequence A001414\cite{oeisA001414}).

\begin{theorem}
If $rank(n)=2$, then $n=2^a3^b5^c$ for some $a$, $b$, $c$ with $a+b+c>0$ and $c<6$. 
\end{theorem}

\begin{proof}[Jānis Iraids, Rihards Opmanis]
If $n$ is divisible by a prime $p>5$, then (since $\f{p}<p$) one can convert the expression of Lemma \ref{primedec} into a shorter expression representing $n$. Hence, if $rank(n)=2$, then $n=2^a3^b5^c$ for some $a, b, c$ with $a+b+c>0$.
But $\f{5^6} < 5\cdot 6$ with the shortest expression being of height 5, as seen in Figure \ref{fig:rank}. Hence, if $n$ is divisible by $5^6$, then, again, one can convert the expression of Lemma \ref{primedec} into a shorter expression representing $n$. 
\end{proof}

\begin{hypo}
\label{genhypo}
All and only numbers of rank $2$ are the numbers $2^a3^b5^c > 5$ with $c<6$.
An equivalent form, extending Hypothesis \ref{powtwohypo}:
\[\f{2^a3^b5^c}=2a+3b+5c\]
for all $a, b, c$ with $a+b+c>0$ and $c<6$ (moreover, the product of $1+...+1$s is shorter than any other representation of $2^a3^b5^c$).
\end{hypo}

As verified by Jānis Iraids, $\f{2^a3^b5^c}=2a+3b+5c$ is true for all $2^a3^b5^c\leq 10^{12}$ with $a+b+c>0$ and $c<6$.

For $c=0$, Hypothesis \ref{genhypo} appears in \cite{raw}.


\subsection{Logarithmic complexity}
\label {subsec:logcomp}
\pdfbookmark[2]{2.3. Logarithmic complexity}{Logarithmic complexity}

Because of Theorem \ref{cbounds}, the values of $\frac{\f{n}}{\log_3{n}}$ are located within the segment $[3, 4.755]$. Hence, the

\begin{definition}
The \textbf{logarithmic complexity} of $n>1$ is defined as
\[\fl{n}=\frac{\f{n}}{\log_3{n}}.\]
\end{definition}

For example, $\fl{2}=\frac{\f{2}}{\log_3{2}} \approx 3.1699$.

\begin{hypo}
\label{genhypo1}
It would follow from Hypothesis \ref{genhypo}, that the values of $\fl{n}$ are \textbf{dense} across the segment $[3, \fl{2}]$, i.e. no subsegment of it is free from the values of $\fl{n}$.
\end{hypo}

\begin{proof}
According to Hypothesis \ref{genhypo}, $\f{2^a3^b}=2a+3b$, for all $a+b>0$. Hence, by choosing $a$ and $b$ appropriately, one can locate the value of $\f{2^a3^b}_{\log}$ within any given subsegment of $[3, \fl{2}]$.
\end{proof}

As observed by Mārtiņš Opmanis (and confirmed by Jānis Iraids for all the numbers $\leq 10^{12}$), it seems, the largest values of $\fl{n}$ are taken by single numbers, see Table \ref{logc}. The lists in braces represent Cunningham chains of primes \cite{chains}.

\begin{center}
\begin{longtable}{|r|c|c|c|l|}
\caption{Largest values of $\fl{n}$}
\label{logc} \\
\hline
$n$ & $\f{n}$ & $\approx \fl{n}$ & $rank(n)$ &Other properties \\
\hline
\endfirsthead
$1439$ & $26$ & $3.928$ & $9$ & $e(26), r(9), \{89, 179, 359, 719, 1439, 2879\}$ \\ 
$23$ & $11$ & $3.854$ & $5$ & $e(11), r(5), \{2, 5, 11, 23, 47\}$ \\ 
$719$ & $23$ & $ 3.841$ & $7$ & $e(23)$, \{89, 179, 359, 719, 1439, 2879\} \\
$179$ & $18$ & $3.812$ & $7$ & $e(18), r(7)$, \{89, 179, 359, 719, 1439, 2879\} \\
$4283$ & $29$ & $3.809$ & $7$ & $e(29), \{2141, 4283\}$ \\
$1438$ & $25$ & $3.777$ & $8$ & $e(25)$, $2×719$ \\
$59$ & $14$ & $3.772$ & $5$ & $e(14), \{29, 59\}$ \\
$6299$ & $30$ & $3.767$ & $7$ & $e(30)$, prime \\
$15287$ & $33$ & $3.763$ & $9$ & $e(33), \{3821, 7643, 15287\}$ \\
$107$ & $16$ & $3.762$ & $5$ & $e(16), \{53, 107\}$ \\
$347$ & $20$ & $3.756$ & $7$ & $e(20), \{173, 347\}$ \\
$1499$ & $25$ & $3.756$ & $7$ & prime\\
$467$ & $21$ & $3.754$ & $5$ & $e(21), \{233, 467\}$ \\
$11807$ & $32$ & $3.749$ & $7$ & $e(32), \{5903, 11807\}$ \\
$263$ & $19$ & $3.746$ & $5$ & $e(19), \{131, 263\}$ \\
$21599$ & $34$ & $3.743$ & $7$ & $e(34), \{2699,5399,10799,21599\}$ \\ \hline
\end{longtable}
\end{center}

The values of $\fl{n}$ become (and, it seems, stay) less than $3.60$ approximately at $n=2\cdot 10^{9}$.

Let's consider the subsegments $[C, D], C<D$ of the segment $[3, 4.755]$ that do not contain the values of $\fl{n}$ at all. Of course, according to Hypothesis \ref{genhypo1}, then $C > \fl{2}$. Let's denote by $C_1$ the infimum of these numbers. I.e. $C_1$ is the point separating the area where the values of $\fl{n}$ are dense, from the area, where these values are not dense. 

\begin{hypo}
\label{hypoc1}
It would follow from Hypothesis \ref{genhypo}, that $C_1 \geq \fl{2}$.
\end{hypo}

On the other hand, for some numbers $C$, $\fl{n}>C$ only for finitely many values of $n$. Let's denote by $C_2$ the infimum of these numbers. This is also known as $\limsup\limits_{n \rightarrow \infty} \fl{n}$. 
I.e. $C_2$ is the point separating the area where the values of $\fl{n}$ are "absolutely sparse", from the area, where these values are not sparse. 

Of course, $C_1 \leq C_2$. Hence,

\begin{hypo}
\label{hypoc2}
It would follow from Hypothesis \ref{genhypo}, that $C_2 \geq \fl{2}$.
\end{hypo}
 More about the possible value of $C_2$ --  in Section \ref{subsec:eexp}.

\section{Analytical results}
\pdfbookmark[1]{3. Analytical results}{Analytical results}
\subsection{Complexity of $2^n-1$}
\pdfbookmark[2]{3.1. Complexity of 2\^{}n-1}{Complexity of $2^n-1$}

For the sake of brevity let us introduce $A(n)=\f{2^n-1}-2n$ and $B(n)=\f{2^n+1}-2n$. We can then establish the following facts.

\begin{theorem}
For $n\geq 1$,
\begin{enumerate}[a)]
\item $A(2n) \leq A(n)+B(n)$;
\item $A(3n) \leq A(n)+B(n)+1$;
\item $A(n+1) \leq A(n)+1$.
\end{enumerate}
\end{theorem}
\begin{proof}
We shall provide the expressions that will result in these upper bounds.
\begin{enumerate}[a)]
\item $2^{2n}-1=(2^n-1)(2^n+1)$. If we take the complexity of both sides and subtract $4n$ we get:
\[A(2n)=\f{2^{2n}-1}-4n \leq \f{2^n-1}-2n+\f{2^n+1}-2n = A(n)+B(n).\]
\item $2^{3n}-1=(2^n-1)((2^n+1)2^n+1)$. Similarly, we get
\[\begin{split} A(3n)&=\f{2^{3n}-1}-6n \leq \\ &\leq \f{2^n-1}-2n+\f{2^n+1}-2n + \f{2^n} - 2n + 1 \leq \\ &\leq A(n)+B(n) + 1. \end{split}\]
\item $2^{n+1}-1=(2^n-1)\cdot 2+1$. Once again we have
\[A(n+1)=\f{2^{n+1}-1}-2n-2 \leq \f{2^n-1}-2n+2-2+1 = A(n)+1.\]
\end{enumerate}
This method can be extended for numbers other than 2, but then it yields significantly less interesting results because of the very inefficient ``$n+1$'' step.
\end{proof}

\begin{corollary}
\label{kaspars}
 (Kaspars Balodis) If $n > 1$, then
\[\f{ 2^n - 1 } \leq 2n + \left\lfloor \log_2{n} \right\rfloor + H(n) - 3,\]
where $H(n)$ is the number of $1$-s in the binary representation of $n$, i.e., $H(n)$ is the Hamming weight of $n$ in binary.
\end{corollary}

\begin{proof}
Using the above theorem we can obtain an upper bound of $A(n)$, setting $B(n)=1$ and $A(1)=A(2)=-1$ and ignoring the rule b) altogether. We will use rule c) per every $1$ in the binary representation of $n$ except the leftmost digit, in total, $H(n)-1$ times. We will use rule a) per every digit, except the two leftmost digits, in total, $\left\lfloor \log_2{n} \right\rfloor -1$ times. In this way we will reduce $n$ to 2 at which point $A(n)=-1$ having ``paid'' 1 for each application of any of the rules a) and c). To sum up:
\[A(n)\leq H(n)-1+\left\lfloor \log_2{n} \right\rfloor -1+(-1).\]
\end{proof}

\begin{corollary}
If $n > 1$ then $\f{ 2^n - 1 } \leq 2n + 2 \left\lfloor \log_2{n}  \right\rfloor - 2$.
\end{corollary}

\subsection{Connection of rank and defect}
\pdfbookmark[2]{3.2. Connection of rank and defect}{Connection of rank and defect}

\begin{definition}
(Harry Altman, Joshua Zelinsky) The \textbf{defect of a number} $n$ is
\[d(n)\overset{\underset{\mathrm{def}}{}}{=} \f{n}-3\cdot \log_3{n}. \]
\end{definition}

\begin{proposition} (Jānis Iraids)
\[d(n) \geq \left\lfloor \frac{rank(n)-1}{2} \right\rfloor \left(1+3\log_3{\frac{6}{7}}\right). \]
\end{proposition}

\begin{proof}
We will prove this by induction on the rank. First of all, for all $n$ having $rank(n)<3$ the proposition is true trivially since
\[\left\lfloor \frac{rank(n)-1}{2} \right\rfloor \left(1+3\log_3{\frac{6}{7}}\right) \leq 0.\]
Now assuming that it is true for all $n$ having $rank(n)<r$. Suppose $r$ is an even number, then we again trivially have
\[\left\lfloor \frac{r-1}{2} \right\rfloor \left(1+3\log_3{\frac{6}{7}}\right) = \left\lfloor \frac{r-2}{2} \right\rfloor \left(1+3\log_3{\frac{6}{7}}\right). \]
If on the other hand, $r$ is an odd number greater than $2$, then the shortest expression for $n$ has height at least $3$, addition being the outermost operation. Now write $n$ as a sum of numbers of even rank lower than $r$ and order them in non-increasing fashion. This can be done since the numbers of odd rank can be split down further merging them with the $n$'s outermost addition.
\[n = \sum_{i=1}^{k}{a_i}, \quad rank(a_1)\geq rank(a_2) \geq \cdots \geq rank(a_k). \] 
Note that $rank(a_1)=r-1$ and $k\geq 2$. For defects, we have
\[d(n) = \sum_{i=1}^{k}{d(a_i)}+3\log_3{\frac{\prod_{i=1}^k{a_i}}{\sum_{i=1}^k{a_i}}}. \]

Following the induction the defect of $a_1$ is at least $\frac{r-3}{2} \left(1+3\log_3{\frac{6}{7}}\right)$. If $rank(a_k)\geq 2$ and so all $a_i$ must necessarily be at least $6$. The expression $3\log_3{\frac{\prod_{i=1}^k{a_i}}{\sum_{i=1}^k{a_i}}}$ is minimised when $k=2$ and $a_1=a_2=6$ producing a minimum of $3$. However, if there are $l$ numbers of rank 0: $rank(a_{k-l+1})=\ldots =rank(a_k)=0$ then $d(a_{k-l+1})=\ldots =d(a_k)=1$ and $\sum_{i=k-l+1}^k{d(a_i)}+3\log_3{\frac{\prod_{i=1}^k{a_i}}{\sum_{i=1}^k{a_i}}}$ is minimised at $l=1$, $a_1=6$ and $a_2=1$ giving the $1+3\log_3{\frac{6}{7}}$. Consequently,
\[d(n) \geq \frac{r-3}{2} \left(1+3\log_3{\frac{6}{7}}\right) + \min{\left(3,1+3\log_3{\frac{6}{7}}\right)} = \]
\[ = \frac{r-1}{2} \left(1+3\log_3{\frac{6}{7}}\right). \]
\end{proof}

\section{Algorithms for computing the complexity of a number}
\pdfbookmark[1]{4. Algorithms for computing the complexity of a number}{Algorithms for computing the complexity of a number}

The purpose of this chapter is to describe several algorithms for computing $\f{n}$. This is useful for both exploring the behaviour of the complexity and as well as a tool for verifying whether some hypotheses about the value of $\f{n}$ hold for large $n$.

It appears that computing the complexity of number $n$ is a relatively difficult task. At the time of writing we are unaware of any algorithm that provably works faster than $O(n^{\log_2{3}-1})$. Yet there exists an algorithm, that in practice runs in approximately $\Theta \left(n \log{n} \right)$ time.

The \textbf{first and simplest algorithm} is essentially as follows: evaluate all the possible expressions
of $\left\{1,+, \cdot\right\}$ with increasing number of ones until an expression giving $n$ is found. We can do so taking as basis the postfix notation, also known as Reverse Polish notation. Note, that in the postfix notation the first (leftmost) symbol will always be 1.

The possible expression generator will use three main rules:
\begin{enumerate}[a)]
{\setlength\itemindent{25pt} \item After 1 only + or 1 can follow. $\cdot $ can not follow because that would mean multiplication with 1 that is meaningless;}
{\setlength\itemindent{25pt} \item After + only $\cdot $ or 1 can follow. + can not follow because that would mean that it was possible to make the previous addition earlier and thus it is not postfix notation;}
{\setlength\itemindent{25pt} \item After $\cdot$ only + or 1 can follow. $\cdot$ can not follow because that would mean that it was possible to make the previous multiplication earlier and thus it is not postfix notation.}
\end{enumerate}

For each total number of ones $x$ in the expression excluding the leading one we can count the number of distinct expressions as $4^x$ -- postfix notation of an expression with $x+1$ ones contains $2x+1$ symbols and for each symbol there are two possible symbols that can follow. Assuming, that evaluation of the expression takes roughly $x$ time, the overall time complexity can be estimated as $O\left(4^x\right)$. Thus the overall running time of the algorithm will be $O\left(4^{\f{n}}\right)$, which taking into account $\f{n} \leq 3 \log_2{n}$ will yield $O\left(n^6\right)$.

However, the expression generating algorithm can be improved by eliminating similar expressions of the same number. For instance, $6$ can be represented as
\begin{itemize}
  \item $(1 + 1) \cdot (1 + 1 + 1)$ -- postfix notation $11+11+1+\cdot$ or
  \item $(1 + 1 + 1) \cdot (1 + 1)$ -- postfix notation $11+1+11+\cdot$.
\end{itemize}
Since multiplication is associative, we can omit one of these representations restricting that the first multiplier has to be smaller or equal to the second multiplier. The same principle applies to additions.

The \textbf{second algorithm} uses the idea of sieving. We will use an array $f[1..n]$ of integer values where the complexity of the number $n$ will be stored as $f[n]$. Before the main routine starts, precalculate the values of $E(n)$ in an array $E[1..m]$.

\algnotext{EndFor}
\algnotext{EndIf}
\algnotext{EndWhile}
\begin{algorithmic}
\For {$i=1$ to $n$}
  \State $f[i] \gets i$
\EndFor
\State $updated\gets True$
\State $height \gets 2$
\While {$updated$}
  \State $updated \gets False$
  \If {$height \equiv 1 \pmod 2$}
    \For {$i=2$ to $n$}
      \State {$a \gets \frac{n-\sqrt{n^2-4E[f[n]]}}{2}$}
      \For {$j=1$ to $a$}
        \If {$f[i] > f[j]+f[i-j]$}
          \State {$f[i] \gets f[j]+f[i-j]$}
          \State {$a \gets \frac{n-\sqrt{n^2-4E[f[n]]}}{2}$}
          \State {$updated \gets True$}
        \EndIf
      \EndFor
    \EndFor
  \Else
    \For {$i=2$ to $\lfloor\sqrt{n}\rfloor$}
      \State {$j \gets i+i$}
      \State {$k \gets 2$}
      \While {$j \leq n$}
        \If {$f[j] > f[k]+f[i]$}
          \State {$f[j] \gets f[k]+f[i]$}
          \State {$updated \gets True$}
        \EndIf
        \State $j \gets j+i$
        \State $k \gets k+1$
      \EndWhile
    \EndFor
  \EndIf
  \State $height \gets height + 1$
\EndWhile
\end{algorithmic}

If the array $f$ is initialized as shown each pass will only update the complexity of numbers that have $rank(n)\geq height$. As a consequence, this algorithm produces the rank of a number; the rank can be stored if necessary. Note that the array $f$ could be initialized with the upper bound provided by Theorem \ref{cbounds} because at no point would $f[i]$ be smaller than $\f{i}$ and yet exceed the value at the corresponding point in unmodified algorithm. To further reduce the running time, one can use a bootstrapping step where numbers that could potentially be used as the smallest of two addends are computed. These are exactly the numbers that cannot be represented best as sums.

While the second algorithm is as fast as any we know, it uses a linear amount of memory -- the array where the complexity values are stored. For $n > 10^{11}$ the calculation thus becomes unfeasible. We used this algorithm to calculate $\f{n}$ and $rank(n)$ up to $n=1.5\cdot 10^9$.  

The \textbf{third and final algorithm} -- the one we used for calculating $\f{n}$ for n up to $10^{12}$ -- is conceptually very simple: for every natural number compute the complexity by definition and store for subsequent steps.
\begin{align}
\f{1}&=1\\
\f{n}&= \min_{a+b=n \vee a\cdot b=n} \{ \f{a}+\f{b}\}
\end{align}
The techniques used are identical to what Fuller \cite{fuller} describes in the comments of his program. For factorisation we used an approach similar to the one described in \cite{pgen}. The core idea is to maintain a priority queue of the so called eliminators -- at any point $n$ for each prime the priority is the smallest integer multiple that is no less than $n$.

\section{Experimental results}
\pdfbookmark[1]{5. Experimental results}{Experimental results}

\subsection{$e(n)$ -- the least number of complexity $n$}
\label{subsec:eexp}
\pdfbookmark[2]{5.1. e(n) - the least number of complexity n}{e(n) - the least number of complexity n}

Function $e(n)$ corresponds to the sequence A005520 \cite{oeisA005520}.
Our observations up $n=89$ are represented in Table \ref{table:e}.

\begin{hypo}
\label{eprime1}
$e(n)$ is prime for all $n$, except $n\in \{1, 4, 7, 11, 25\}$. \end{hypo}

\begin{obs}
\label{eprime2}
For $k \leq 3$, the number $\frac{e(n)-k}{k+1}$ is prime for almost all $n$.
\end{obs}

If, for $k=1$, Observation \ref{eprime2} holds for an infinite number of values of $n$, it would imply that there is an infinite number of \textbf{Sophie Germain primes} -- these are defined as integers $p$ such that $p$ and $2p+1$ are both primes.

Moreover, it seems that the sequence of {e(n)} contains primes which are the end numbers of increasingly long Cunningham chains. \textbf{Cunningham chain} (CC) of length $k$ is defined as a sequence of $k$ primes $\{p_1,p_2,...,p_k\}$ such that $p_{i+1}=2\cdot p_i+1, 1\leq i < k$ \cite{chains}. 
In particular,
\begin{itemize}
\item $e(13)$ is the end number of the first CC of length 5: $\{2, 5, 11, 23, 47\}$;
\item $e(26)$ is the end number of another CC of length 5: $\{89, 179, 359, 719, 1439\}$;
\item $e(27)$ is the end number of the first CC of length 6: 
  \[\{89, 179, 359, 719, 1439, 2879\};\]
\item $e(80)$ is the end number of another CC of length 5.
\end{itemize}

The above-mentioned are the only CCs of length $\geq 5$ backward-generated by $e(n), n\leq 89$.

For the following $19$ values of $n\leq 89$, $e(n)$ generates CCs of length 4:

\[ \{11, 23, 34, 49, 51, 60, 61, 65, 66, 67, 70, 72, 73, 74, 77, 84, 86, 87, 89\}.\]

If, for $k=2; 3$, Observation \ref{eprime2} holds for an infinite number of values of $n$, it would imply that there is an infinite number of integers $p$ such that $p$ and $3p+2$ (or, correspondingly, $p$ and $4p+3$) are both primes \cite{dickson}.

\begin{center}
\begin{longtable}{|c|r|r|r|l|}
\caption{Prime factorizations of numbers close to $e(n)$}
\label{table:e} \\
\hline
\multicolumn{1}{|c|}{$n$} & \multicolumn{1}{c|}{$\frac{e(n)-2}{3}$} & \multicolumn{1}{c|}{$\frac{e(n)-1}{2}$} & \multicolumn{1}{c|}{$e(n)$} & \multicolumn{1}{c|}{$e(n)+1$} \\ \hline
\endfirsthead

1 & -- & -- & $1$ & $2$\\
2 & -- & -- & $2$ & $3$\\
3 & -- & $1$ & $3$ & $2^2$\\
4 & -- & -- & $2^2$ & $5$\\
5 & $1$ & $2$ & $5$ & $2\cdot3$\\
6 & $5$ & $3$ & $7$ & $2^3$\\
7 & -- & -- & $2\cdot5$ & $11$\\
8 & $3$ & $5$ & $11$ & $2^2\cdot3$\\
9 & $5$ & $2^3$ & $17$ & $2\cdot3^2$\\
10 & -- & -- & $2\cdot11$ & $23$\\
11 & $7$ & $11$ & $23$ & $2^3\cdot3$\\
12 & $13$ & $2^2\cdot5$ & $41$ & $2\cdot3\cdot7$\\
13 & $3\cdot5$ & $23$ & $47$ & $2^4\cdot3$\\
14 & $19$ & $29$ & $59$ & $2^2\cdot3\cdot5$\\
15 & $29$ & $2^2\cdot11$ & $89$ & $2\cdot3^2\cdot5$\\
16 & $5\cdot7$ & $53$ & $107$ & $2^2\cdot3^3$\\
17 & $5\cdot11$ & $83$ & $167$ & $2^3\cdot3\cdot7$\\
18 & $59$ & $89$ & $179$ & $2^2\cdot3^2\cdot5$\\
19 & $3\cdot29$ & $131$ & $263$ & $2^3\cdot3\cdot11$\\
20 & $5\cdot23$ & $173$ & $347$ & $2^2\cdot3\cdot29$\\
21 & $5\cdot31$ & $233$ & $467$ & $2^2\cdot3^2\cdot13$\\
22 & $227$ & $11\cdot31$ & $683$ & $2^2\cdot3^2\cdot19$\\
23 & $239$ & $359$ & $719$ & $2^4\cdot3^2\cdot5$\\
24 & $11\cdot37$ & $13\cdot47$ & $1223$ & $2^3\cdot3^2\cdot17$\\
25 & -- & -- & $2\cdot719$ & $1439$\\
26 & $479$ & $719$ & $1439$ & $2^5\cdot3^2\cdot5$\\
27 & $7\cdot137$ & $1439$ & $2879$ & $2^6\cdot3^2\cdot5$\\
28 & $5\cdot251$ & $7\cdot269$ & $3767$ & $2^3\cdot3\cdot157$\\
29 & $1427$ & $2141$ & $4283$ & $2^2\cdot3^2\cdot7\cdot17$\\
30 & $2099$ & $47\cdot67$ & $6299$ & $2^2\cdot3^2\cdot5^2\cdot7$\\
31 & $3359$ & $5039$ & $10079$ & $2^5\cdot3^2\cdot5\cdot7$\\
32 & $5\cdot787$ & $5903$ & $11807$ & $2^5\cdot3^2\cdot41$\\
33 & $5\cdot1019$ & $7643$ & $15287$ & $2^3\cdot3\cdot7^2\cdot13$\\
34 & $23\cdot313$ & $10799$ & $21599$ & $2^5\cdot3^3\cdot5^2$\\
35 & $3\cdot3733$ & $107\cdot157$ & $33599$ & $2^6\cdot3\cdot5^2\cdot7$\\
36 & $5\cdot23\cdot131$ & $2\cdot11299$ & $45197$ & $2\cdot3^6\cdot31$\\
37 & $18679$ & $28019$ & $56039$ & $2^3\cdot3\cdot5\cdot467$\\
38 & $5\cdot5443$ & $40823$ & $81647$ & $2^4\cdot3^6\cdot7$\\
39 & $32999$ & $49499$ & $98999$ & $2^3\cdot3^2\cdot5^3\cdot11$\\
40 & $54419$ & $81629$ & $163259$ & $2^2\cdot3^2\cdot5\cdot907$\\
41 & $53\cdot1283$ & $101999$ & $203999$ & $2^5\cdot3\cdot5^3\cdot17$\\
42 & $80627$ & $120941$ & $241883$ & $2^2\cdot3^2\cdot6719$\\
43 & $5\cdot24763$ & $185723$ & $371447$ & $2^3\cdot3^2\cdot7\cdot11\cdot67$\\
44 & $180179$ & $270269$ & $540539$ & $2^2\cdot3^3\cdot5\cdot7\cdot11\cdot13$\\
45 & $196799$ & $295199$ & $590399$ & $2^6\cdot3^2\cdot5^2\cdot41$\\
46 & $302399$ & $453599$ & $907199$ & $2^6\cdot3^4\cdot5^2\cdot7$\\
47 & $173\cdot2083$ & $540539$ & $\num{1081079}$ & $2^3\cdot3^3\cdot5\cdot7\cdot11\cdot13$\\
48 & $617039$ & $925559$ & $\num{1851119}$ & $2^4\cdot3^3\cdot5\cdot857$\\
49 & $680399$ & $\num{1020599}$ & $\num{2041199}$ & $2^4\cdot3^6\cdot5^2\cdot7$\\
50 & $\num{1081079}$ & $\num{1621619}$ & $\num{3243239}$ & $2^3\cdot3^4\cdot5\cdot7\cdot11\cdot13$\\
51 & $\num{1280159}$ & $\num{1920239}$ & $\num{3840479}$ & $2^5\cdot3^3\cdot5\cdot7\cdot127$\\
52 & $\num{2187359}$ & $\num{3281039}$ & $\num{6562079}$ & $2^5\cdot3^3\cdot5\cdot7^2\cdot31$\\
53 & $\num{2735519}$ & $\num{4103279}$ & $\num{8206559}$ & $2^5\cdot3^2\cdot5\cdot41\cdot139$\\
54 & $\num{3898919}$ & $\num{5848379}$ & $\num{11696759}$ & $2^3\cdot3^2\cdot5\cdot32491$\\
55 & $\num{4882919}$ & $\num{7324379}$ & $\num{14648759}$ & $2^3\cdot3^2\cdot5\cdot7\cdot5813$\\
56 & $13\cdot59\cdot9697$ & $\num{11156399}$ & $\num{22312799}$ & $2^5\cdot3^3\cdot5^2\cdot1033$\\
57 & $\num{9164959}$ & $\num{13747439}$ & $\num{27494879}$ & $2^5\cdot3\cdot5\cdot7^3\cdot167$\\
58 & $\num{13915439}$ & $\num{20873159}$ & $\num{41746319}$ & $2^4\cdot3^3\cdot5\cdot7\cdot11\cdot251$\\
59 & $\num{17417399}$ & $\num{26126099}$ & $\num{52252199}$ & $2^3\cdot3^2\cdot5^2\cdot7\cdot11\cdot13\cdot29$\\
60 & $\num{26110559}$ & $\num{39165839}$ & $\num{78331679}$ & $2^5\cdot3^2\cdot5\cdot7\cdot19\cdot409$\\
61 & $\num{36202319}$ & $\num{54303479}$ & $\num{108606959}$ & $2^4\cdot3^3\cdot5\cdot7\cdot11\cdot653$\\
62 & $6577\cdot7247$ & $\num{71495279}$ & $\num{142990559}$ & $2^5\cdot3^2\cdot5\cdot109\cdot911$\\
63 & $\num{67699439}$ & $\num{101549159}$ & $\num{203098319}$ & $2^4\cdot3^3\cdot5\cdot17\cdot5531$\\
64 & $\num{91328639}$ & $\num{136992959}$ & $\num{273985919}$ & $2^7\cdot3^2\cdot5\cdot13\cdot3659$\\
65 & $\num{127340639}$ & $\num{191010959}$ & $\num{382021919}$ & $2^5\cdot3^4\cdot5\cdot7\cdot4211$\\
66 & $\num{165145679}$ & $\num{247718519}$ & $\num{495437039}$ & $2^4\cdot3^3\cdot5\cdot7^2\cdot31\cdot151$\\
67 & $\num{227109119}$ & $\num{340663679}$ & $\num{681327359}$ & $2^8\cdot3^2\cdot5\cdot7^2\cdot17\cdot71$\\
68 & $\num{335430119}$ & $\num{503145179}$ & $\num{1006290359}$ & $2^3\cdot3^2\cdot5\cdot601\cdot4651$\\
69 & $\num{468798119}$ & $\num{703197179}$ & $\num{1406394359}$ & $2^3\cdot3^3\cdot5\cdot7\cdot17\cdot31\cdot353$\\
70 & $\num{619264799}$ & $\num{928897199}$ & $\num{1857794399}$ & $2^5\cdot3^3\cdot5^2\cdot7\cdot11\cdot1117$\\
71 & $\num{909474719}$ & $\num{1364212079}$ & $\num{2728424159}$ & $2^5\cdot3^2\cdot5\cdot7\cdot11^2\cdot2237$\\
72 & $\num{1247732639}$ & $\num{1871598959}$ & $\num{3743197919}$ & $2^5\cdot3^4\cdot5\cdot7\cdot11^3\cdot31$\\
73 & $40499\cdot41221$ & $\num{2504113919}$ & $\num{5008227839}$ & $2^9\cdot3^3\cdot5\cdot7\cdot11\cdot941$\\
74 & $\num{2290896719}$ & $\num{3436345079}$ & $\num{6872690159}$ & $2^4\cdot3^3\cdot5\cdot7\cdot454543$\\
75 & $\num{3279830399}$ & $7643\cdot643693$ & $\num{9839491199}$ & $2^7\cdot3^4\cdot5^2\cdot7\cdot11\cdot17\cdot29$\\
76 & $\num{4495159679}$ & $\num{6742739519}$ & $\num{13485479039}$ & $2^7\cdot3^2\cdot5\cdot11^3\cdot1759$\\
77 & $\num{5574925439}$ & $\num{8362388159}$ & $\num{16724776319}$ & $2^7\cdot3^8\cdot5\cdot7\cdot569$\\
78 & $\num{8226486239}$ & $\num{12339729359}$ & $\num{24679458719}$ & $2^5\cdot3^2\cdot5\cdot7\cdot19\cdot128861$\\
79 & $\num{11841566159}$ & $\num{17762349239}$ & $\num{35524698479}$ & $2^4\cdot3^2\cdot5\cdot49339859$\\
80 & $\num{14737208639}$ & $\num{22105812959}$ & $\num{44211625919}$ & $2^6\cdot3^2\cdot5\cdot7^2\cdot11\cdot19\cdot1499$\\
81 & $\num{20797230719}$ & $\num{31195846079}$ & $\num{62391692159}$ & $2^7\cdot3^2\cdot5\cdot7^2\cdot221059$\\
82 & $\num{31251071039}$ & $\num{46876606559}$ & $\num{93753213119}$ & $2^6\cdot3^2\cdot5\cdot7^2\cdot61\cdot10891$\\
83 & $\num{40517305919}$ & $\num{60775958879}$ & $\num{121551917759}$ & $2^6\cdot3^4\cdot5\cdot7\cdot13\cdot29\cdot1777$\\
84 & $\num{54513320399}$ & $\num{81769980599}$ & $\num{163539961199}$ & $2^4\cdot3^3\cdot5^2\cdot7\cdot11\cdot196657$\\
85 & $\num{83528413919}$ & $\num{125292620879}$ & $\num{250585241759}$ & $2^5\cdot3^4\cdot5\cdot7\cdot2762183$\\
86 & $\num{106809776639}$ & $\num{160214664959}$ & $\num{320429329919}$ & $2^9\cdot3^3\cdot5\cdot7^2\cdot37\cdot2557$\\
87 & $\num{141615840239}$ & $\num{212423760359}$ & $\num{424847520719}$ & $2^4\cdot3^4\cdot5\cdot7\cdot9366127$\\
88 & $\num{210123688319}$ & $\num{315185532479}$ & $\num{630371064959}$ & $2^7\cdot3^3\cdot5\cdot7\cdot13\cdot17\cdot23581$\\
89 & $\num{290857880879}$ & $\num{436286821319}$ & $\num{872573642639}$ & $2^4\cdot3^3\cdot5\cdot7^3\cdot19\cdot61987$\\ \hline 
\end{longtable}
\end{center} 

The behaviour of $e(n)$ provides some evidence that the logarithmic complexity of n does not tend to 3, and even that the constant $C_2=\limsup\limits_{n \rightarrow \infty} \fl{n}$ is greater than $\fl{2}$. It is useful to note the following fact:

\begin{proposition}
\[C_2=\limsup\limits_{n \rightarrow \infty} \fl{n}=\limsup\limits_{n \rightarrow \infty} \fl{e(n)}\]
\end{proposition}
\begin{proof}
Since $\fl{e(n)}$ is a subsequence of $\fl{n}$ the $\geq$ follows obviously.

The $\leq $ is proven from the contrary. Assume that there exists an infinite list of numbers $\{x_i\}$, that only have a finite number of elements from sequence $e(n)$ of logarithmic complexity $\geq \min_i \{\fl{x_i}\}$. Furthermore, there is an infinite subsequence $\{y_i\}$ of $\{x_i\}$'s such that the complexity of numbers strongly increases. But any number of this sequence has a corresponding number in $e(n)$, namely $e(\f{y_i})$ that has greater or equal logarithmic complexity.

Alternatively, on the axis of logarithmic complexity, $e(n)$ yields the rightmost point of numbers of complexity $n$.
\end{proof}

\begin{hypo}
\label{enlog1}
The limit $\lim\limits_{n \rightarrow \infty} \fl{e(n)}$ exists and
\[\lim\limits_{n \rightarrow \infty} \fl{e(n)} = \limsup\limits_{n \rightarrow \infty} \fl{e(n)}\]
\end{hypo}

\begin{figure}
\centering
\resizebox{0.95\textwidth}{!}{\input{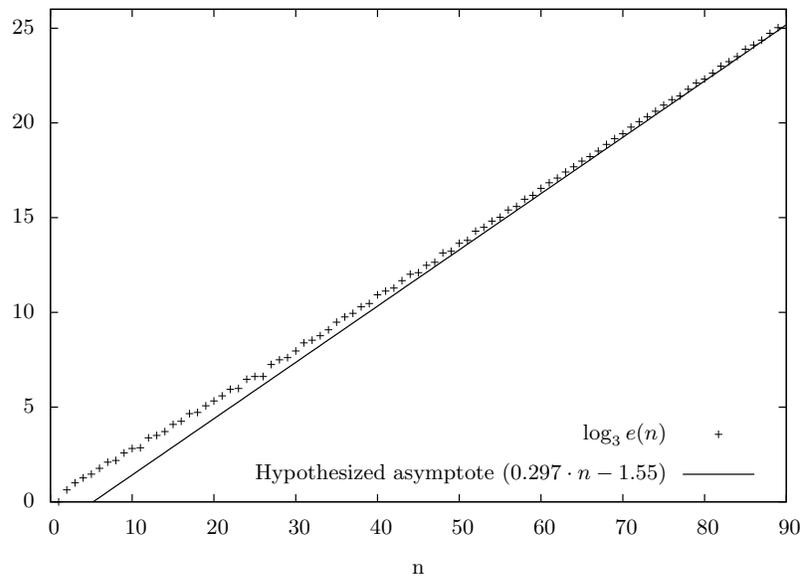}}
\caption{Values of $\log_3 e(n)$ compared with hypothesized asymptote}
\label{fig:easymptote}
\end{figure}

\begin{figure}
\centering
\resizebox{0.95\textwidth}{!}{\input{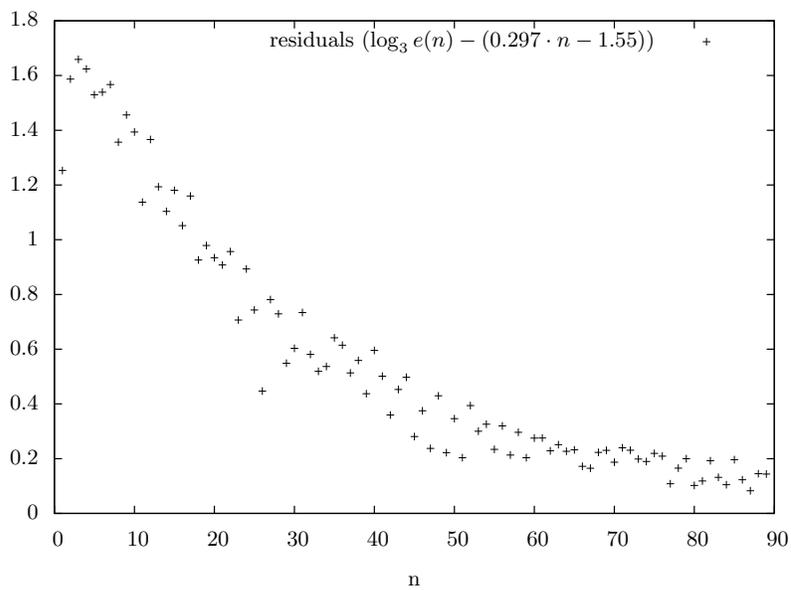}}
\caption{Residuals of the asymptote}
\label{fig:eresiduals}
\end{figure}

As one can see in Figures \ref{fig:easymptote} and \ref{fig:eresiduals}: 

\begin{hypo}
\label{enlog2}
$\log_3{e(n)}$ behaves almost linearly, namely,
\[\log_3{e(n)} \approx 0.297n - 1.55.\]
\end{hypo}

Hence, 

\[\fl{e(n)} = \frac{n}{\log_3{e(n)}} \approx 3.37 + \frac{5.2}{0.297 n-1.55},\]

and $\lim\limits_{n \rightarrow \infty} \fl{e(n)} \approx 3.37$.

\begin{lemma}
For all $n>1$:

\[\fl{n} \leq \fl{e(\f{n})}.\]
\end{lemma}

\begin{proof}
Obviously, $\f{n}=\f{e(\f{n})}$, and $n \geq e(\f{n})$. Hence,
$\log_3{n} \geq \log_3{e(\f{n})}$, and $\fl{n} \leq \fl{e(\f{n})}$.
\end{proof} 

Thus, if $\epsilon>0$, then, it seems, $\fl{n} > 3.37+\epsilon$ can be true only for finitely many values of $n$. Thus, in terms of Section \ref{subsec:logcomp}:

\begin{hypo}
\label{hypoc2prim}
It would follow from Hypothesis \ref{enlog2}, that $C_2 \leq 3.37$ (approximately). 
\end{hypo}

\subsection{Structure of shortest expressions}
\pdfbookmark[2]{5.2. Structure of shortest expressions}{Structure of shortest expressions}

Analyzing the structure of the shortest expressions representing numbers, we have come to the conclusion that we should not only abandon any attempts to obtain a "deterministic" method for construction of shortest expressions, and turn to "nondeterministic" methods. We should abandon also radical attempts to predict regularities in the structure of shortest expressions. 

For example, one might propose the following hypothesis: 
if $p$ is the smallest prime divisor of $n>1$, then
\[\f{n}=\min{(1+\f{n-1}, \f{p}+\f{n/p})}.\]

For $p=n$, and $p=2$ (with prime $n/p$) this hypothesis appears in Guy \cite{guy}: for any prime $p$, 
\[\f{p}=1+\f{p-1},\] 
\[\f{2p}=\min{(1+\f{2p-1}, 2+\f{p})}.\]
Similarly, one could suppose: 
\[\f{3p}=\min{(1+\f{3p-1}, 3+\f{p})}.\]

The first hypothesis fails, the smallest counterexample being

\[p=\num{353942783} = 2 \cdot 3 +(1+2^2 \cdot 3^2)(2+3^4(1+2 \cdot 3^{10})),\]
$\f{p}=63, 1+\f{p-1}=64$, found by Martin N. Fuller, 2008, see \cite{oeisA005245}.

The second hypothesis fails as well, the smallest counterexample being
\[2p=\num{10278600694} = 2 \cdot 3 +(1+2 \cdot 3^2(1+2^4))(1+3^{14}(1+2 \cdot 3)),\]
$\f{2p}=72, 1+\f{2p-1}=2+\f{p}=73$, found by Jānis Iraids, 2010. 

By analogy, it seems, the third hypothesis also should fail, but this does not happen for $3p\leq 10^{12}$. 

Thus, when trying to build the shortest expression representing a number, subtraction of 1 and division by primes are not universal candidates for the first operation. How about subtraction of other numbers?

Subtractions of 3, 4, 5, 7, 10, 11, 13, etc. include subtraction of 1. Thus, it remains to consider only subtractions of 6, 8, 9, 12, etc.

The first number, for which subtraction of 6 is necessary as the first operation, is the above prime found by Fuller:

\[\num{353942783} = 2 \cdot 3 +(1+2^2 \cdot 3^2)(2+3^4(1+2 \cdot 3^{10})).\]

Until $10^{12}$, there are only 21360 numbers for which subtraction of 6, 8 or 9 is necessary as the first operation.

Until $10^{12}$, there are exactly 3 numbers for which the first operation must be subtraction of 8:
\[\num{341317451698}=2 \cdot prime=2^3+(1+2^4(1+2 \cdot 3^3))(1+3^{18});\]
\[\num{474934483834}=2 \cdot 6011 \cdot 39505447=2^3+(1+2^4)(1+2^3 \cdot 3^4)(1+3^{16});\]
\[\num{782747233558}=2 \cdot prime=2^3+(1+2^4(1+2 \cdot 3^5))(1+3^{15}(1+2 \cdot 3)).\]

Until $10^{12}$, there are 119 numbers for which the first operation must be subtraction of 9, the first three ones being:
\[\num{16534727299}=103 \cdot 160531333=3^2+(1+2^7 \cdot 3^3)(1+3^{14});\]
\[\num{68238632999}=prime=3^2+(1+1+3^2)(1+2^4 \cdot 3^4)(1+3^{14});\]
\[\num{85619928299}=prime=3^2+(1+2^2 \cdot 3)(1+2^4)(1+3^{18}).\]

Necessity for subtraction of 12 (or larger addendum) was not detected for numbers until $10^{12}$.

According to Corollary \ref{smallestaddend}, if $n\geq 29$ and the shortest expression for $n$ is a sum $n=a+b$, then the smaller addendum $a\leq 2n^{\log_2{3}-1} \approx 2n^{0.585}$. However, the above observations show that for $n\leq 10^{12}$ the smaller addendum does not exceed 9.

\subsection{Complexity of $2^n+1$ and $2^n-1$}
\pdfbookmark[2]{5.3. Complexity of 2\^{}n-1}{Complexity of $2^n-1$ and $2^n-1$}

Since, it seems, $\f{2^n}=2n$, one can suppose

\begin{hypo}
\label{hypo2ap1}
For all $n \geq 0$, except 3 and 9, 
\[\f{2^n+1}=2n+1.\]
\end{hypo}
Hypothesis \ref{hypo2ap1} is true for all $2^n+1 \leq 10^{12}$, i.e. for all $n \leq 39$ – as verified by Jānis Iraids.

Both exceptions are due to relatively massive divisibility by 3:
\[2^3+1 = 9 = 2 \cdot 2 \cdot 2+1=3 \cdot 3;\] 
\[2^9+1 = 513 = (3 \cdot 3 \cdot 2+1) \cdot 3 \cdot 3 \cdot 3. \]

On the other hand, since we do not have subtraction in our expression basis $\{1,+,\cdot \}$,  the numbers $2^n-1$ seem to be more complicated than $2^n$. In Theorem \ref{kaspars}, an upper bound of  $\f{2^n-1}$ was proved. In Table \ref{kaspars2} this result is compared with experimental data for $n \leq 39$.

\begin{longtable}{|c|c|c|c|c|c|}
\caption{Complexity of $\f{2^n-1}$}
\label{kaspars2} \\
\hline
$n$ & $\f{2^n-1}-2n$ & $\lfloor \log_2{n}\rfloor + H(n) - 3$ & $n$ & $\f{2^n-1}-2n$ & $\lfloor \log_2{n}\rfloor + H(n) - 3$ \\ 
\hline
\endfirsthead
$1$ & $-1$ & $-$ & $21$ & $2$ & $4$ \\
$2$ & $-1$ & $-1$ & $22$ & $3$ & $4$ \\
$3$ & $0$ & $0$ & $23$ & $4$ & $5$ \\
$4$ & $0$ & $0$ & $24$ & $2$ & $3$ \\
$5$ & $1$ & $1$ & $25$ & $3$ & $4$ \\
$6$ & $0$ & $1$ & $26$ & $3$ & $4$ \\
$7$ & $1$ & $2$ & $27$ & $2$ & $5$ \\
$8$ & $1$ & $1$ & $28$ & $3$ & $4$ \\
$9$ & $1$ & $2$ & $29$ & $4$ & $5$ \\
$10$ & $2$ & $2$ & $30$ & $3$ & $5$ \\
$11$ & $3$ & $3$ & $31$ & $4$ & $6$ \\
$12$ & $1$ & $2$ & $32$ & $3$ & $3$ \\
$13$ & $2$ & $3$ & $33$ & $4$ & $4$ \\
$14$ & $2$ & $3$ & $34$ & $4$ & $4$ \\
$15$ & $2$ & $4$ & $35$ & $4$ & $5$ \\
$16$ & $2$ & $2$ & $36$ & $2$ & $4$ \\
$17$ & $3$ & $3$ & $37$ & $2$ & $5$ \\
$18$ & $1$ & $3$ & $38$ & $2$ & $5$ \\
$19$ & $2$ & $4$ & $39$ & $3$ & $6$ \\
$20$ & $2$ & $3$ & $-$ & $-$ & $-$ \\ \hline
\end{longtable}

Thus, it seems, the upper bound of Theorem \ref{kaspars} is exact for all $n=2^k, k>0$:

\begin{obs}
For all $0 \leq k \leq 5$,  
\[\f{2^{2^k}-1}=2 \cdot 2^k + k - 2.\]
\end{obs}

\subsection{Observing ranks}
\pdfbookmark[2]{5.4. Observing ranks}{Observing ranks}
\label{subsec:exprank}

The first number for which subtraction of 6 is necessary as the first operation to obtain the shortest expression of minimum height, is
\[\num{22697747}=prime=2\cdot 3 +(2+3^7)(1+2^7\cdot 3^4),\]
complexity 55, rank 5.

The values of $r(n)$ up to $n=19$ are represented in Table \ref{table:rn}. The lists in braces represent Cunningham chains of primes \cite{chains}.

\begin{longtable}{|c|l|c|l|}
\caption{$r(n)$ -- the least number of rank $n$}
\label{table:rn} \\
\hline
$n$ & $r(n)$ & $\f{r(n)}$ & Other properties \\ 
\hline
\endfirsthead
$1$ & $2$ & $2$ & $e(2), \{2, 5, 11, 23, 47\}$ \\
$2$ & $6$ & $5$ & $2\cdot3$ \\
$3$ & $7$ & $6$ & $e(6), \{3, 7\}$ \\
$4$ & $14$ & $8$ & $2\cdot7$ \\
$5$ & $23$ & $11$ & $e(11), \{2, 5, 11, 23, 47\}$ \\
$6$ & $86$ & $14$ & $2\cdot43$  \\
$7$ & $179$ & $18$ & $e(18), \{89, 179, 359, 719, 1439, 2879\}$ \\
$8$ & $538$ & $21$ & $2\cdot269$ \\
$9$ & $1439$ & $26$ & $e(26), \{89, 179, 359, 719, 1439, 2879\}$ \\
$10$ & $9566$ & $30$ & $2\cdot4783$ \\ 
$11$ & $21383$ & $33$ & $\{10691, 21383, 42767\}$ \\
$12$ & $122847$ & $37$ & $3\cdot40949$ \\
$13$ & $777419$ & $44$ & prime \\
$14$ & $1965374$ & $46$ & $2\cdot982687$ \\
$15$ & $6803099$ & $51$ & $\{3401549, 6803099\}$ \\
$16$ & $19860614$ & $53$ & $2\cdot9930307$ \\
$17$ & $26489579$ & $55$ & $\{13244789, 26489579, 52979159\}$  \\
$18$ & $269998838$ & $61$ & $2\cdot4093\cdot32983$\\
$19$ & $477028439$ & $64$ & $14207\cdot33577$ \\
\hline
\end{longtable}

\begin{figure}
\centering
\resizebox{1.0\textwidth}{!}{\input{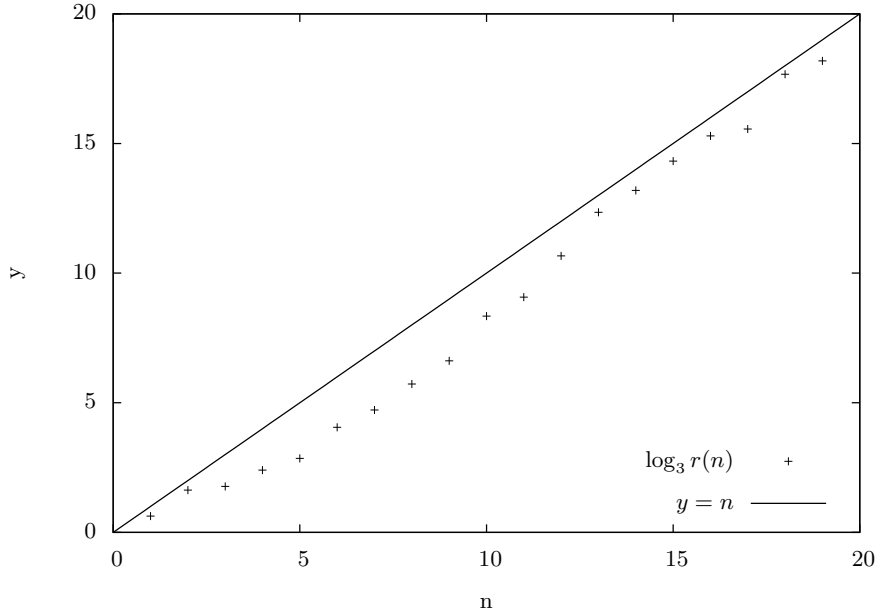}}
\caption{Values of $\log_3{r(n)}$ compared with $n$.
}
\label{fig:rasymptote}
\end{figure}

In Figure \ref{fig:rasymptote}, the values of $\log_3{r(n)}$ are compared with $n$.

\begin{obs}
$\log_3{r(n)}$ tends to $n$, hence, it seems, $r(n) \approx 3^n$.
\end{obs} 

\subsection{Collapse of powers}
\pdfbookmark[2]{5.5. Collapse of powers}{Collapse of powers}
\label{subsec:collapse}

While attempting to prove or disprove the Hypothesis \ref{powtwohypo} one might try to generalize the hypothesis by looking for other numbers $n$, such that the shortest expressions for all the integer powers of $n$: $n^2$, $n^3$, $n^4$, $\ldots$ can be obtained as products of shortest expressions for $n$.

\begin{definition}
If this property holds, then the number $n$ is called \textbf{resistant}.
\end{definition}

Obviously, if $n$ is resistant, then for all $k>0$,

\[\f{n^k}=k \f{n}; \]
\[\fl{n^k}=\fl{n}. \]

Currently, only the number $3$ is proved to be resistant. According to Hypothesis \ref{powtwohypo}, another resistance candidate is the number $2$. Existence of resistant numbers other than powers of $3$ would provide a lower bound on $C_2$:

\begin{proposition}
If $n=2$ or $n \neq 3^k$ is resistant, $C_2 \geq \fl{n} > 3$.
\end{proposition}

\begin{definition}
If $k>0$ is the least number such that $\f{n^k} < k \f{n}$, let us say that the number $n$ \textbf{collapses at $k$}.
\end{definition}

For example, the number $5$ collapses at $6$ (see Section \ref{subsec:ranking}):
\[\f{5^6}=29 < 6\f{5}=6\cdot 5 = 30. \]

However, it seems, most primes collapse already at 2, for example,
\[\f{11^2}=15 < 2\f{11}=2\cdot 8 = 16; 11^2=121=1+2^3\cdot 3\cdot 5. \]

Of the $168$ primes until $10^3$, $120$ primes collapse at $2$ ($71\%$), $24$ -- at $3$ ($14\%$), $3$ -- at $4$ ($2\%$), the remaining $21$ do not collapse until $4$ ($13\%$), and, currently, 12 of them are not known to collapse at all, namely,
\[2, 3, 109, 163, 379, 433, 487, 541, 577, 733, 739, 811. \]
For details, see Table \ref{table:collapse1000}. The rest of numbers all collapse at $2$ or $3$, and their logarithmic complexity exceeds $3.417$. It seems, primes having larger logarithmic complexity are more likely to collapse.

\begin{longtable}{|c|c|c|c|c|}
\caption{Powers in which primes until 1000 collapse, sorted by increasing logarithmic complexity}
\label{table:collapse1000} \\
\hline
$p$ & $collapses$ & $\f{p}$ & $rank(p)$ & $\approx\fl{p}$ \\ 
\hline
\endfirsthead
$3$ & $-$ & $3$ & $1$ & $3.000$ \\
$2$ & $?(>39)$ & $2$ & $1$ & $3.170$ \\
$487$ & $?(>4)$ & $18$ & $3$ & $3.196$ \\
$163$ & $?(>5)$ & $15$ & $3$ & $3.235$ \\
$433$ & $?(>4)$ & $18$ & $3$ & $3.257$ \\
$109$ & $?(>5)$ & $14$ & $3$ & $3.278$ \\
$811$ & $?(>4)$ & $20$ & $3$ & $3.280$ \\
$577$ & $?(>4)$ & $19$ & $3$ & $3.283$ \\
$769$ & $3$ & $20$ & $3$ & $3.307$ \\
$757$ & $\leq 6$ & $20$ & $5$ & $3.314$ \\
$541$ & $?(>4)$ & $19$ & $3$ & $3.317$ \\
$739$ & $?(>4)$ & $20$ & $5$ & $3.326$ \\
$73$ & $6$ & $13$ & $3$ & $3.329$ \\
$379$ & $?(>4)$ & $18$ & $5$ & $3.331$ \\
$733$ & $?(>4)$ & $20$ & $5$ & $3.331$ \\
$271$ & $4$ & $17$ & $3$ & $3.334$ \\
$193$ & $4$ & $16$ & $3$ & $3.340$ \\
$991$ & $\leq 12$ & $21$ & $5$ & $3.344$ \\
$37$ & $5$ & $11$ & $3$ & $3.347$ \\
$977$ & $2$ & $21$ & $5$ & $3.351$ \\
$19$ & $6$ & $9$ & $3$ & $3.358$ \\
$97$ & $6$ & $14$ & $3$ & $3.362$ \\
$257$ & $\leq 6$ & $17$ & $3$ & $3.366$ \\
$937$ & $3$ & $21$ & $5$ & $3.372$ \\
$673$ & $3$ & $20$ & $5$ & $3.374$ \\
$919$ & $3$ & $21$ & $5$ & $3.381$ \\
$181$ & $3$ & $16$ & $3$ & $3.381$ \\
$661$ & $3$ & $20$ & $5$ & $3.384$ \\
$7$ & $9$ & $6$ & $3$ & $3.387$ \\
$653$ & $2$ & $20$ & $5$ & $3.390$ \\
$337$ & $3$ & $18$ & $5$ & $3.398$ \\
$641$ & $4$ & $20$ & $3$ & $3.400$ \\
$883$ & $3$ & $21$ & $5$ & $3.401$ \\
$127$ & $2$ & $15$ & $5$ & $3.402$ \\
$881$ & $2$ & $21$ & $5$ & $3.402$ \\
$877$ & $3$ & $21$ & $5$ & $3.405$ \\
$241$ & $3$ & $17$ & $3$ & $3.405$ \\
$631$ & $2$ & $20$ & $5$ & $3.408$ \\
$457$ & $3$ & $19$ & $5$ & $3.408$ \\
$331$ & $3$ & $18$ & $5$ & $3.408$ \\
$5$ & $6$ & $5$ & $1$ & $3.413$ \\
\hline
\end{longtable}

Of the $1229$ primes until $10^4$, $1030$ primes collapse at $2$ ($84\%$), $122$ -- at $3$ ($10\%$), the remaining $77$ do not collapse until $3$ ($6\%$).  

Of the $78498$ primes until $10^6$, $71391$ primes collapse at $2$ ($91\%$), the remaining $7107$ do not collapse until $2$ ($9\%$), and, currently, only 2722 are not known to collapse at all ($4\%$).

The following observations are true for primes less than $10^6$.
\begin{obs}
Almost all primes collapse at $2$.
\end{obs}

\begin{obs}
If a prime $p$ does not collapse at all, then $\fl{p} < 3.364$. If a prime $p$ collapses, then $\fl{p} > 3.180 > \fl{2}$.
\end{obs}

If we turn to composite numbers, we can encounter numbers that collapse because one ore more of their prime divisors collapse -- even when their initial shortest expression did not directly contain the prime as a multiplier. $3^4+1=82$ is an example of such a number:
\begin{itemize}
\item $\f{82}=13$ and $\f{41}=12$
\item $\f{82^{12}}<12\cdot 13$, because $\f{41^{12}}\leq 131 <12\cdot 11$:
\[\begin{split} 41^{12}=&[2^3\cdot 3^2\cdot \left(2^2+1\right)\left(3^7+1\right) \\
&\left(2\cdot 3^4\cdot \left(2^3\cdot 3^2\cdot \left(2\cdot 3^3+1\right)+1\right)+1\right) \\
& \left(2\cdot 3^4\cdot \left(2\cdot 3^9\cdot \left(2\cdot 3+1\right)+1\right)+1\right)]+1
\end{split}\]
\end{itemize}
On the other hand, there could possibly be composite numbers that do not collapse even though some of their prime divisors do. Obviously, a necessary condition is that the shortest expression for the composite number does not directly contain the collapsing prime as a multiplier. One candidate is $3^5+1=244=2^2\cdot 61$; $\f{244}=16$ but $\f{61}=13$. Using a heuristic algorithm Juris Čerņenoks was able to produce expressions for $61^k$ for $k=7..15$. His results suggest that $61$ does not collapse well enough, i.e., the number of ones saved from collapsing is less than required to catch up with the expression $244=3^5+1$.

On a side note, the above mentioned $82$ and $244$ are interesting because if they proved to be resistant, then $C_2$ would exceed $\fl{2}$.

\section{Conclusions}
\pdfbookmark[1]{6. Conclusions}{Conclusions}

Trying to explore representing of natural numbers by arithmetical expressions using ones, addition, multiplication and parentheses, one arrives very soon at the problems that are easy to formulate, but (it seems) extremely hard to solve. 

Consider, for example, the above Hypothesis \ref{powtwohypo} stating that the best way of representing of $2^n$ is $(1+1)(1+1)...$. We consider proving or disproving of Hypothesis \ref{powtwohypo} as one of the biggest challenges of number theory.

Almost as challenging seems Hypothesis \ref{unlrank} stating that 
as a function of $n$, $rank(n)$ is unlimited. Rank is an additional (to $\f{n}$) measure of integer complexity introduced in this paper.

As another challenge we regard determining of the distribution of the values of logarithmic complexity $\fl{n}$ which are located within the segment $[3, 4.755]$. First, denote by $C_1$ the point separating the area (on the left) where the values of $\fl{n}$ are dense, from the area, where these values are not dense. On the other hand, on the right, the values of $\fl{n}$ are "absolutely sparse": for most $C$, $\fl{n}>C$ only for finitely many values of $n$. Denote by $C_2$ the point separating the area (on the right) where the values of $\fl{n}$ are "absolutely sparse", from the area, where these values are not sparse. Of course, $C_1 \leq C_2$. Our Hypotheses \ref{genhypo} and \ref{enlog2} (if true) imply that

\[3.1699 \approx \fl{2} \leq C_1 \leq C_2 \leq 3.37.\]

Our main experimental "device" is the database containing the values of $\f{n}$ up to $n=10^{12}$ calculated by Jānis Iraids. The database can be accessed by using an online calculator page linked from \cite{oeisA005245}.  

And finally, our Hypothesis \ref{eprime1} and Observation \ref{eprime2} (if true) imply that there is an infinite number of Sophie Germain primes, and even that there is an infinite number of Cunningham chains of length 4 (at least).

\bibliographystyle{splncs}

\phantomsection
\addcontentsline{toc}{chapter}{References}
\bibliography{exp}

\end{document}